# $\alpha$-GAUSS CURVATURE FLOWS

LAMI KIM AND KI-AHM LEE

Abstract. In this paper, we study the deformation of the n-dimensional strictly convex hypersurface in $\mathbb{R}^{n+1}$ whose speed at a point on the hypersurface is proportional to $\alpha$-power of Gauss curvature. For $\frac{1}{n} < \alpha \leq 1$, we prove that there exist the strictly convex smooth solutions if the initial hypersurface is strictly convex and smooth and the solution hypersurfaces converge to a point. We also show the asymptotic behavior of the rescaled hypersurfaces, in other words, the rescaled manifold converges to a strictly convex smooth manifold. Moreover, there exists a subsequence whose limit satisfies a certain equation.

## 1. Introduction

We study the regularity of the $\alpha$-Gauss curvature flow, describing the deformation of an n-dimensional compact strictly convex hypersurface $\Sigma$ in $\mathbb{R}^{n+1}$ which is worn down by collision from any random angle. For example, the stone under the river is eroded by water, where the speed of erosion is in proportion to the $\alpha$-Gauss curvature $K^\alpha$. Let $X(\cdot,\cdot) : \Sigma \times [0,T) \to \mathbb{R}^{n+1}$ be an embedding and set $\Sigma_t = X(\Sigma, t)$. Then the hypersurface evolves by the following flow:

$$\begin{aligned} \frac{\partial X}{\partial t}(x,t) &= -K^\alpha(x,t)\,\nu(x,t) \\ X(x,0) &= X_0(x) \end{aligned} \tag{1.1}$$

where $\nu$ denotes the unit outward normal to $\Sigma_t$ and $\alpha > 0$.

**1.1.** Let $(0, T^*]$ be the maximal interval in which $vol(\Sigma_t)$ is nonzero. Then there are the known results for the case of $\alpha = 1$ in the flow following (1.1). If the initial hypersurface in $\mathbb{R}^3$ is the smooth and strictly convex and has the central symmetry, then the solution $\Sigma_t$ converges to a point as a spherical shape [F]. Also Tso, [T], showed existence and regularity of the solution when the initial hypersurface embedded in $\mathbb{R}^{n+1}$ is the smooth and strictly convex. In other words, the solution $\Sigma_t$ preserves the smoothness and convexity in the time interval $(0, T^*]$. For a smooth, compact, and strictly convex initial hypersurface in $\mathbb{R}^3$, the solution hypersurface $\Sigma_t$ converges to a point and the rescaled solution hypersurface $\tilde{\Sigma}_t$ approaches the round sphere with normalized volume and for a non-smooth initial hypersurface, the viscosity solution has $C^{1,1}$-regularity in the time interval $(0, T^*)$ and $C^\infty$-regularity for $t \geq t_0$ where $t_0$ depends on the volume and diameter of the initial hypersurface $\Sigma_0$ [A1].

For $\alpha = \frac{1}{n+2}$, the solution, $\Sigma_t$, is known as an affine normal flow. There exists a unique, smooth, convex solution such that the hypersurfaces $\Sigma_t$ converge to a point and the rescaled solution converges to an ellipsoid if the initial hypersurface is a compact, smooth, and strictly convex [A5]. For $\frac{1}{n+2} < \alpha \leq \frac{1}{n}$ or $0 < \alpha \leq \frac{1}{n}$ under the assumption that the isoperimetric ratios are bounded, there exist a smooth, strictly convex solution converging to a point and a rescaled solution satisfying a certain equation [A2]. In addition, for $\alpha = \frac{1}{n}$, the rescaled solution converges to a sphere and this holds for $\alpha \geq \frac{1}{n}$ if the initial hypersurface is very close to a sphere [C1].

**1.2. Main Theorem.** Let us denote the rescaled $\Sigma$ and a support function $S$ by $\tilde{\Sigma}$ and $\tilde{S}$ respectively so that the enclosed volume becomes normalized. We state the main theorem.

**Theorem 1.1.**
*Let $\Sigma_0 = X(\Sigma, 0)$ be a compact, connected, strictly convex smooth manifold in $\mathbb{R}^{n+1}$. Assume $\frac{1}{n} < \alpha \leq 1$. Then*

*(i) there exist a time $T^*$ and strictly convex smooth solutions $\{\Sigma_t = X(\Sigma, t)\}$ satisfying* (1.1) *for $t \in [0, T^*)$, and $\Sigma_t$ converges to a point as $t$ approaches to $T^*$.*

*(ii) The principal curvatures of the rescaled hypersurfaces $\tilde{\Sigma}$ have the uniform upper and lower bounds. In other words, let us denote the eigenvalues of $(\tilde{h}^i_j)$ by $\tilde{\lambda}_k$ for $k = 1, \ldots, n$ and the smallest and largest one by $\tilde{\lambda}_{min}$ and $\tilde{\lambda}_{max}$, respectively. Then we have*

$$\frac{1}{M} \leq \tilde{\lambda}_{min} \leq \tilde{\lambda}_{max} \leq M$$

*for some constant $0 < M < \infty$.*

*(iii) For any sequence $\tau_i \to \infty$, there exists a subsequence $\tau_{i_k}$ such that the rescaled manifold $\tilde{\Sigma}_{\tau_{i_k}}$ converges to a strictly convex manifold $\tilde{\Sigma}_{T^*}$ uniformly in $C^\infty$-norm.*

*(iv) In addition, the limit, $\tilde{S}_*(\cdot)$, of the volume normalized solution $\tilde{S}(\cdot, \tau_{i_k})$ satisfies the equation $\tilde{K}_*^\alpha = \tilde{C}_* \tilde{S}_*$ a.e. for some positive constant $\tilde{C}_*$, where $\tilde{K}_*$ is the Gauss curvature of $\tilde{\Sigma}_{T^*}$.*

**1.3. Outline.** Each sections in this paper will be organized as follows. In this section, we have introduced the known results for Gauss curvature flow and our main theorem. In section 2, we state the definitions of a metric, the second fundamental form, some curvatures and the support function. In addition, we obtain the evolution equations for the geometric quantities. In section 3, we prove the hypersufaces preserve the strict convexity and we also get the uniform upper bounds of the Gauss curvature $K$ and the eigenvalues of the reverse second fundamental form. Then using these bounds, we can get the uniform bound of curvatures of hypersurface $\Sigma$. This is proved in Section 5.2. An integral quantity plays the key role in getting the asymptotic behavior of hypersurface and $C^{1,1}$-regularity of the rescaled solution. Also, the curvature bounds of the rescaled hypersurfaces will be introduced in section 4. In [T] and [C1], the authors showed that there exists the finite time $T^*$ such that the solution $\Sigma_t$ converges to a point as $t$ approaches to $T^*$. In the last section, we shall discuss existence of solutions and the asymptotic behavior of the rescaled hypersurfaces. Throughout the whole section, we consider the case $\frac{1}{n} < \alpha \leq 1$ unless there is some explicit assumption on $\alpha$. We will also assume that $\Sigma_t$ is smooth whenever we prove a priori-estimates.

## 2. Evolution of the metric and curvature

**2.1.** Let $\{x_1, \ldots, x_n\}$ be the local coordinates of $\Sigma_t$ and $\nu$ be the outward unit normal vector to $\Sigma_t$. Then the induced metric and the second fundamental form are defined by

$$g_{ij} = \left\langle \frac{\partial X}{\partial x^i}, \frac{\partial X}{\partial x^j} \right\rangle \quad \text{and} \quad h_{ij} = -\left\langle \frac{\partial^2 X}{\partial x^i \partial x^j}, \nu \right\rangle.$$

Also, the Weingarten map $\mathcal{W}_p : T_pM \to T_pM$ for the hypersurface $M \subset \mathbb{R}^{n+1}$ can be given by

$$h^i_j = g^{ik} h_{kj},$$

where $(g^{ij})$ denotes the inverse matrix of $(g_{ij})$, and then $\sigma_k = \sum_{1 \leq i_1 < \cdots < i_k \leq n} \lambda_{i_1} \lambda_{i_2} \cdots \lambda_{i_k}$, $H = \text{trace}(h^i_j) = \sigma_1 = \sum_{1 \leq i \leq n} \lambda_i$, $K = \det(h^i_j) = \sigma_n = \lambda_1 \lambda_2 \cdots \lambda_n$, and $|A|^2 = h_{ij} h^{ij} = \lambda_1^2 + \cdots + \lambda_n^2$, where $\lambda_1, \ldots, \lambda_n$ are the eigenvalues of the Weingarten map at $p$.

The evolutions of the metric, the second fundamental form, and curvature are the following. Throughout this paper, the symbol $\Box$ will be used in place of the operator $K^\alpha (h^{-1})^{kl} \nabla_k \nabla_l$. For the detailed proof, the readers can refer to Chapter 2, [Z].

**Lemma 2.1.** *Let $\Sigma_0$ be strictly convex and $\Sigma_t = X(\Sigma, t)$ be smooth. For the $\alpha$-Gauss curvature flow, we have*

*(i)* $$\frac{\partial g_{ij}}{\partial t} = -2K^\alpha h_{ij}$$

*(ii)* $$\frac{\partial \nu}{\partial t} = g^{ij} \frac{\partial K^\alpha}{\partial x^i} \frac{\partial X}{\partial x^j} = \nabla^j K^\alpha \frac{\partial X}{\partial x^j}$$

*(iii)* $$\frac{\partial h_{ij}}{\partial t} = \nabla_i \nabla_j K^\alpha - K^\alpha h_{jk} h^k_i$$

$$= \alpha \Box h_{ij} + \alpha^2 K^\alpha (h^{-1})^{kl}(h^{-1})^{mn}\nabla_i h_{kl}\nabla_j h_{mn} - \alpha K^\alpha (h^{-1})^{km}(h^{-1})^{nl}\nabla_i h_{mn}\nabla_j h_{kl}$$
$$+ \alpha K^\alpha H h_{ij} - (1+n\alpha)K^\alpha h_{jl}h_i^l$$

*(iv)* $\dfrac{\partial K}{\partial t} = \alpha \Box K + \alpha(\alpha-1)K^{\alpha-1}(h^{-1})^{ij}\nabla_i K \nabla_j K + K^{\alpha+1}H$

*(v)* $\dfrac{\partial K^\alpha}{\partial t} = \alpha \Box K^\alpha + \alpha K^{2\alpha} H$

*(vi)* $\dfrac{\partial H}{\partial t} = \alpha \Box H + \alpha^2 K^{\alpha-2} g^{ij}\nabla_i K \nabla_j K - \alpha K^\alpha g^{ij}(h^{-1})^{km}(h^{-1})^{nl}\nabla_i h_{mn}\nabla_j h_{kl} + \alpha K^\alpha H^2$
$+(1-n\alpha)K^\alpha |A|^2$

*(vii)* $\dfrac{\partial |X|^2}{\partial t} = \Box |X|^2 - 2K^\alpha (h^{-1})^{kl} g_{kl} + 2(n-1)K^\alpha \langle X, \nu\rangle$

**2.2.** The support function $S(z,t)$ of the strictly convex hypersurface is given by

$$S(z) = \langle z, X(\nu^{-1}(z),t)\rangle, \quad \text{for} \quad z \in \mathbb{S}^n, \tag{2.1}$$

where $\mathbb{S}^n$ denotes a unit sphere. Then $X(z)$ can be written as

$$X(z) = S(z)z + \overline{\nabla} S(z)$$

from the definition of the support function, (2.1), and $\overline{\nabla}_i S(z) = \langle X(z), \overline{\nabla}_i z\rangle$ for the connection of the standard metric $\overline{g}$ on $\mathbb{S}^n$. We also have

$$\frac{\partial z}{\partial x^i} = h_{ik} g^{kl} \frac{\partial X}{\partial x^l} \tag{2.2}$$

from the relation between the tangent vector and the normal vector and the definition of the second fundamental form. In addition,

$$h_{ij} = \overline{\nabla}_i \overline{\nabla}_j S + S\overline{g}_{ij} \tag{2.3}$$

where $\overline{g}_{ij}$ is the metric on $\mathbb{S}^n$, which this can be obtained by taking covariant derivatives of (2.1), [Z].

Now we have the following properties and the evolution equations (see Chapter 3 of [Z] to understand the proof in detail).

**Lemma 2.2.** *Let $\Sigma_0$ be strictly convex and $X(\nu^{-1}(z),t)$ be smooth, where $z \in \mathbb{S}^n$. For the $\alpha$-Gauss curvature flow, we have*

*(i)* $\overline{g}_{ij} = h_{ik} g^{kl} h_{lj}$ *and* $\overline{g}^{ij} = (h^{-1})^{ik} g_{kl} (h^{-1})^{lj}$

*(ii)* $h^i_j = (h^{-1})^{ik}\overline{g}_{kj}$

*(iii)* $H = \overline{g}_{ik}(h^{-1})^{ki}$, $|A|^2 = g^{kl}\overline{g}_{kl}$, *and* $K = \det(h^i_j) = \dfrac{\det(\overline{g}_{ij})}{\det(\overline{\nabla}_i\overline{\nabla}_j S + S\overline{g}_{ij})}$

*(iv) Set* $\mathcal{S}_k$ *be the k-th symmetric polynomial of* $h_{ij}$ *while* $\sigma_k$ *is the k-th symmetric polynomial of* $h^i_j$. *Then* $\mathcal{S}_n = K^{-1}$.

The following lemma gives us the evolution equations of the support function, second fundamental form, and curvatures for the standard metric $\overline{g}_{ij}$ on $S^n$.

**Lemma 2.3.** *Let* $\Sigma_0$ *be strictly convex and* $X(\nu^{-1}(z), t)$ *be smooth, where* $z \in \mathbb{S}^n$. *For the* $\alpha$*-Gauss curvature flow, we have*

*(i)* $\dfrac{\partial S}{\partial t} = -K^\alpha = -\mathcal{K}^{-\alpha}$ *or* $\Big(-\dfrac{\partial S}{\partial t}\Big)\mathcal{K}^\alpha = 1$, *where* $\mathcal{K} = K^{-1}$

*(ii)* $\dfrac{\partial h_{ij}}{\partial t} = -\overline{\nabla}_i\overline{\nabla}_j K^\alpha - K^\alpha\overline{g}_{ij} = -\Big(\overline{\nabla}_i\overline{\nabla}_j\mathcal{K}^{-\alpha} + \mathcal{K}^{-\alpha}\overline{g}_{ij}\Big)$

*(iii)* $\dfrac{\partial H}{\partial t} = g^{ij}\Big(\overline{\nabla}_i\overline{\nabla}_j K^\alpha + K^\alpha\overline{g}_{ij}\Big)$

*(iv)* $\dfrac{\partial |A|^2}{\partial t} = 2\overline{g}_{ij}h^{ij}K^\alpha$

*(v)*
$$\begin{aligned}\frac{\partial K}{\partial t} &= K(h^{-1})^{ij}\Big(\overline{\nabla}_i\overline{\nabla}_j K^\alpha + K^\alpha\overline{g}_{ij}\Big)\\ &= \alpha K^\alpha(h^{-1})^{ij}\overline{\nabla}_i\overline{\nabla}_j K + \alpha(\alpha-1)K^{\alpha-1}(h^{-1})^{ij}\overline{\nabla}_i K\overline{\nabla}_j K + K^{\alpha+1}H\end{aligned}$$

*(vi)* $\dfrac{\partial K^\alpha}{\partial t} = \alpha K^\alpha(h^{-1})^{ij}\overline{\nabla}_i\overline{\nabla}_j K^\alpha + \alpha K^{2\alpha}H$

*(vii)* $\dfrac{\partial \mathcal{K}}{\partial t} = -\mathcal{K}(h^{-1})^{ij}\left(\overline{\nabla}_i\overline{\nabla}_j\mathcal{K}^{-\alpha} + \mathcal{K}^{-\alpha}\overline{g}_{ij}\right)$

*Proof.* Taking the time derivative of (2.1) gives us

$$\frac{\partial S}{\partial t} = \Big\langle z, \nabla X\cdot\frac{\partial \nu^{-1}}{\partial t} + \frac{\partial X}{\partial t}\Big\rangle$$

and then (*i*) comes from (1.1). Also we can obtain (*ii*) and (*iv*) by the definitions of $h_{ij}$ and $|A|^2$, respectively.

We know that

$$\begin{aligned}\frac{\partial}{\partial t}H &= \overline{g}_{ij}(h^{-1})^{ik}(h^{-1})^{lj}(\overline{\nabla}_k\overline{\nabla}_l K^\alpha + \overline{g}_{kl}K^\alpha)\\ &= g^{kl}(\overline{\nabla}_k\overline{\nabla}_l K^\alpha + \overline{g}_{kl}K^\alpha)\end{aligned}$$

by (*ii*). From the evolution equation of the second fundamental form, we get the evolution equation of $K$:

$$\begin{aligned}\frac{\partial K}{\partial t} &= -K\overline{g}^{jm}h_{mi}\overline{g}_{jn}(h^{-1})^{nk}(h^{-1})^{li}\frac{\partial}{\partial t}h^{kl}\\ &= K(h^{-1})^{kl}\overline{\nabla}_k\overline{\nabla}_l K^\alpha + K^{\alpha+1}H,\end{aligned}$$

which implies (*vi*). Also (*vii*) is obtained directly from the definition of $\mathcal{K}$.

□

In addition, $S$ satisfies, as in [T],

$$-S_t(z,t)\Big[\det(\bar{\nabla}_i\bar{\nabla}_j S(z,t) + S(z,t)\delta_{ij})\Big]^\alpha = 1 \quad \text{for } (z,t) \in \mathbb{S}^n \times (0,T^*), \tag{2.4}$$

which comes from Lemma 2.2 (*iii*) and Lemma 2.3 (*i*).

## 3. Curvature Estimate

We define the width, the inner radius and the outer radius of the convex hypersurface as follows:

- the inner radius $r_{in} = \sup\{r : B_r(y) \text{ is enclosed by } X \text{ for some } y \in \mathbb{R}^{n+1}\}$
- the outer radius $r_{out} = \inf\{r : B_r(y) \text{ encloses } X \text{ for some } y \in \mathbb{R}^{n+1}\}$

Now we shall show that the strict convexity of $\Sigma_t$ will be preserved under the flow.

**Lemma 3.1.**
*If $\Sigma_0$ is strictly convex, $\Sigma_t = X(\Sigma, t)$ is also strictly convex for $t > 0$ as long as it is smooth. We also have*

$$\inf_{x\in\Sigma} K(x,t) \geq \inf_{x\in\Sigma} K(x,0) > 0.$$

*Proof.* Let $Z(t) = \inf_{x\in\Sigma} K(x,t)$ and assume that the minimum is achieved at $X = X(x,t)$. Then, at $X$, we have

$$\nabla_i\nabla_j K \geq 0 \quad \text{and} \quad \nabla_i K = 0\,,$$

and hence we get

$$\begin{aligned}\frac{\partial K}{\partial t} &= \alpha \Box K + \alpha(\alpha-1)K^{\alpha-1}(h^{-1})^{ij}\nabla_i K \nabla_j K + K^{\alpha+1}H \\ &\geq K^{\alpha+1}H.\end{aligned} \tag{3.1}$$

Now, $H \geq nK^{1/n}$ implies

$$\frac{\partial Z}{\partial t} \geq nZ^{\alpha+1+1/n}.$$

By the maximum principle, we can get $Z(t) \geq Z(0) > 0$ which gives the positive lower bound of $K$ for $t > 0$, and then the strict convexity of $\Sigma_t$.

□

We have the following lemma (cf. [A2]). We shall use the idea of Lemma 3.5 in [Z].

**Lemma 3.2.**
*Let $\Sigma_0$ be convex, $\Sigma_t = X(\mathbb{S}^n, t)$ be smooth for $t$ in $[0, T^*)$, and $\alpha > 0$. Also let us consider the sphere with radius $r_{in}(T^* - \delta)$ and center at the origin contained in $\Sigma_{T^*-\delta}$ and set $\rho_0 = \frac{1}{2}r_{in}(T^* - \delta)$ where $\delta$ is any positive constant satisfying $\delta < T^*$. Then there is a constant $C > 0$ such that*

$$\sup_{z\in\mathbb{S}^n,\ 0\leq t\leq T^*-\delta} K^\alpha(z,t) \leq C = \max\left(\sup_{z\in\mathbb{S}^n} K^\alpha(z,0), \left(\frac{n\alpha+1}{n\alpha\rho_0}\right)^{n\alpha}\right).$$

*Proof.* We consider the function $\varphi = \dfrac{K^\alpha}{S-\rho_0}$, where $S$ is the support function. Here $S(z,t) = \langle z, X(\nu^{-1}(z), t)\rangle$ and then

$$\frac{\partial S}{\partial t} = \left(z, \frac{\partial X}{\partial t}\right) = \left(z, -K^\alpha \nu\right) = -K^\alpha.$$

Let us assume that $\varphi$ has its maximum at $(z_0,\ t_0)$ for $t_0 \leq T^* - \delta$. Then, at $(z_0,\ t_0)$, we get

$$\varphi_t \geq 0,\ \overline{\nabla}_i\varphi = 0 \text{ and } \overline{\nabla}_i\overline{\nabla}_j\varphi \leq 0.$$

Now we have $0 = \overline{\nabla}_i\varphi = \frac{(S-\rho_0)\overline{\nabla}_i K^\alpha - K^\alpha\overline{\nabla}_i S}{(S-\rho_0)^2} = \frac{\overline{\nabla}_i K^\alpha}{S-\rho_0} - \frac{K^\alpha\overline{\nabla}_i S}{(S-\rho_0)^2}$, so $\overline{\nabla}_i K^\alpha = \frac{K^\alpha\overline{\nabla}_i S}{S-\rho_0}$. Since

$$\begin{aligned}
0 \geq \overline{\nabla}_i\overline{\nabla}_j\varphi &= \overline{\nabla}_i\Big(\frac{\overline{\nabla}_j K^\alpha}{S-\rho_0} - \frac{K^\alpha\overline{\nabla}_j S}{(S-\rho_0)^2}\Big)\\
&= \frac{\overline{\nabla}_i\overline{\nabla}_j K^\alpha}{S-\rho_0} - \frac{\overline{\nabla}_j K^\alpha\overline{\nabla}_i S}{(S-\rho_0)^2} - \frac{\overline{\nabla}_i K^\alpha\overline{\nabla}_j S + K^\alpha\overline{\nabla}_i\overline{\nabla}_j S}{(S-\rho_0)^2} + \frac{2K^\alpha\overline{\nabla}_j S\overline{\nabla}_i S}{(S-\rho_0)^3}\\
&= \frac{\overline{\nabla}_i\overline{\nabla}_j K^\alpha}{S-\rho_0} - \frac{K^\alpha\overline{\nabla}_i\overline{\nabla}_j S}{(S-\rho_0)^2},
\end{aligned}$$

we also have $\overline{\nabla}_i\overline{\nabla}_j K^\alpha \leq \frac{K^\alpha\overline{\nabla}_i\overline{\nabla}_j S}{S-\rho_0}$. Therefore $\varphi$ satisfies, at $(z_0,\ t_0)$,

$$\begin{aligned}
0 \leq \frac{\partial}{\partial t}\varphi &= \frac{1}{S-\rho_0}\Big(\alpha K^\alpha(h^{-1})^{ij}\overline{\nabla}_i\overline{\nabla}_j K^\alpha + \alpha K^{2\alpha}H + \frac{K^{2\alpha}}{S-\rho_0}\Big)\\
&\leq \frac{1}{S-\rho_0}\Big\{\alpha K^\alpha(h^{-1})^{ij}\Big(\frac{K^\alpha\overline{\nabla}_i\overline{\nabla}_j S}{S-\rho_0}\Big) + \alpha K^{2\alpha}H + \frac{K^{2\alpha}}{S-\rho_0}\Big\}.
\end{aligned}$$

From Lemma 2.2, we can derive

$$\begin{aligned}
0 &\leq \frac{\alpha K^\alpha(h^{-1})^{ij}K^\alpha}{S-\rho_0}(h_{ij} - S\overline{g}_{ij}) + \alpha K^{2\alpha}H + \frac{K^{2\alpha}}{S-\rho_0}\\
&= \frac{K^{2\alpha}}{S-\rho_0}\Big(n\alpha - \alpha\rho_0 H + 1\Big),
\end{aligned}$$

which means

$$0 \leq (n\alpha + 1) - \alpha\rho_0 H.$$

If $\frac{n\alpha+1}{\alpha\rho_0} < H$, that is, $H > \frac{C}{\rho_0} > \frac{C}{r_{in}} > \frac{C}{r_{out}}$, where $C = \frac{n\alpha+1}{\alpha}$, we get a contradiction. Therefore $H$ is bounded, so $K^\alpha$ is bounded since $K^\alpha \leq n^{-n\alpha}H^{n\alpha}$. Now we conclude that

$$\sup_{z\in\mathbb{S}^n,\ 0\leq t\leq T^*-\delta} K^\alpha(z,t) \leq C = \max\left(\sup_{z\in\mathbb{S}^n} K^\alpha(z,0), \Big(\frac{n\alpha+1}{n\alpha\rho_0}\Big)^{n\alpha}\right).$$

□

Now we consider the eigenvalues of the reverse second fundamental form.

**Lemma 3.3.**
*Let $\Sigma_0$ be strictly convex, $\Sigma_t = X(\Sigma, t)$ be smooth, and $\frac{1}{n} < \alpha \leq 1$. Also set $\mathcal{H} = (h^{-1})^{ij}g_{ij}$. Then there is a constant $C > 0$ such that*

$$\sup_{x\in\Sigma}\mathcal{H}(x,t) \leq C = \max\left(\sup_{x\in\Sigma}\Big(\frac{n\alpha-1}{\alpha}K^{-1/n}(x,t)\Big), \sup_{x\in\Sigma}\mathcal{H}(x,0)\right)$$

*for $t > 0$ as long as it is smooth.*

*Proof.* First we have the evolution equation for $\mathcal{H}$:

$$\begin{aligned}\frac{\partial \mathcal{H}}{\partial t} &= \alpha\Box\mathcal{H} - 2\alpha K^{\alpha}(h^{-1})^{\gamma\beta}(h^{-1})^{ip}(h^{-1})^{kq}(h^{-1})^{lj}g_{ij}\nabla_{\gamma}h_{kl}\nabla_{\beta}h_{pq} \\ &\quad - \alpha^2 K^{\alpha}(h^{-1})^{ik}(h^{-1})^{lj}(h^{-1})^{mn}(h^{-1})^{pq}g_{ij}\nabla_k h_{mn}\nabla_l h_{pq} \\ &\quad + \alpha K^{\alpha}(h^{-1})^{ik}(h^{-1})^{lj}(h^{-1})^{mp}(h^{-1})^{nq}g_{ij}\nabla_k h_{mn}\nabla_l h_{pq} - \alpha K^{\alpha}H\mathcal{H} + n(1+n\alpha)K^{\alpha} - 2nK^{\alpha}\end{aligned}$$

since we can obtain

$$\alpha\Box\mathcal{H} = -\alpha(h^{-1})^{ik}(h^{-1})^{lj}g_{ij}\Box h_{kl} + 2\alpha K^{\alpha}(h^{-1})^{\gamma\beta}(h^{-1})^{ip}(h^{-1})^{kq}(h^{-1})^{lj}g_{ij}\nabla_{\gamma}h_{kl}\nabla_{\beta}h_{pq} \tag{3.2}$$

from the second derivatives of $\mathcal{H}$ and we also have the Codazzi identity and symmetry of $h_{ij}$. Then at a maximum point we get

$$\begin{aligned}\frac{\partial \mathcal{H}}{\partial t} &\le -2\alpha K^{\alpha}(h^{-1})^{kl}(h^{-1})^{ip}(h^{-1})^{mq}(h^{-1})^{nj}g_{ij}\nabla_k h_{mn}\nabla_l h_{pq} \\ &\quad + \alpha K^{\alpha}(h^{-1})^{ik}(h^{-1})^{lj}(h^{-1})^{mp}(h^{-1})^{nq}g_{ij}\nabla_k h_{mn}\nabla_l h_{pq} \\ &\quad - \alpha K^{\alpha}H\mathcal{H} + n(n\alpha-1)K^{\alpha} \\ &\le -\alpha K^{\alpha}(h^{-1})^{kl}(h^{-1})^{ip}(h^{-1})^{mq}(h^{-1})^{nj}g_{ij}\nabla_k h_{mn}\nabla_l h_{pq} \\ &\quad - \alpha K^{\alpha}H\mathcal{H} + n(n\alpha-1)K^{\alpha} \\ &\le -(\alpha H\mathcal{H} - n(n\alpha-1))K^{\alpha} \\ &\le -(\alpha n K^{1/n}\mathcal{H} - n(n\alpha-1))K^{\alpha}\end{aligned}$$

since $H \ge n(K)^{1/n} \ge c_0 > 0$ for some positive constant $c_0$. On the other hand, we have a contradiction if $\mathcal{H} > \frac{n\alpha-1}{\alpha K^{1/n}}$ at a maximum point. Hence $\mathcal{H} \le \frac{n\alpha-1}{\alpha}K^{-1/n}$. Then the result follows.

□

## 4. Integral Quantities and Asymptotic Behavior of Hypersurface

We shall define the volume $V(t)$ and the area $\mathcal{A}(t)$ enclosed by convex hypersurface $\Sigma$ as follows:

- the volume function $V(t) = \frac{1}{n+1}\int_{\Sigma}\langle X, \nu\rangle\, d\sigma_{\Sigma} = \frac{1}{n+1}\int_{\mathbb{S}^n}\frac{S}{K}\, d\sigma_{\mathbb{S}^n}$
- the area function $\mathcal{A}(t) = \int_{\Sigma} d\sigma_{\Sigma} = \int_{\mathbb{S}^n}\frac{1}{K}\, d\sigma_{\mathbb{S}^n}$

**Lemma 4.1.** *For the strictly convex and smooth solution $\Sigma_t = X(\mathbb{S}^n, t)$ of $\alpha$-Gauss curvature flow* (1.1)*, we have*

$$\frac{\partial}{\partial t} V(t) = -\int_{\mathbb{S}^n} \frac{1}{K^{1-\alpha}}\, d\sigma_{\mathbb{S}^n}\,.$$

*Proof.* First observe that from Lemma 2.2 and Lemma 2.3, we have

$$\begin{aligned}\int_{\mathbb{S}^n} S\mathcal{K}_t\, d\sigma_{\mathbb{S}^n} &= \int_{\mathbb{S}^n} S\mathcal{K}(h^{-1})^{ij}(\overline{\nabla}_i\overline{\nabla}_j S_t + S_t\overline{g}_{ij})\, d\sigma_{\mathbb{S}^n} \\ &= \int_{\mathbb{S}^n} S_t\mathcal{K}(h^{-1})^{ij}(\overline{\nabla}_i\overline{\nabla}_j S + S\overline{g}_{ij})\, d\sigma_{\mathbb{S}^n} \\ &= \int_{\mathbb{S}^n} S_t\mathcal{K}(h^{-1})^{ij} h_{ij}\, d\sigma_{\mathbb{S}^n} = n\int_{\mathbb{S}^n} S_t\mathcal{K}\, d\sigma_{\mathbb{S}^n}\end{aligned}$$

since $\overline{\nabla}_i \mathcal{K}(h^{-1})^{ij} = 0$. Hence we have

$$\begin{aligned}\frac{\partial}{\partial t} V(t) &= \frac{1}{n+1}\int_{\mathbb{S}^n} (\mathcal{K}S_t + S\mathcal{K}_t)\, d\sigma_{\mathbb{S}^n} \\ &= \int_{\mathbb{S}^n} \mathcal{K}S_t\, d\sigma_{\mathbb{S}^n} = -\int_{\mathbb{S}^n} \frac{1}{K^{1-\alpha}}\, d\sigma_{\mathbb{S}^n}.\end{aligned}$$

□

Now let us consider the rescaled solution

$$\tilde{X}(\tau) := \frac{X(t)}{V(t)^{1/(n+1)}} \tag{4.1}$$

and also assume that the normalized volume $\tilde{V}(\tau) = \frac{1}{n+1}\left(\int_{\mathbb{S}^n} \frac{\tilde{S}}{\tilde{K}}\, d\sigma_{\mathbb{S}^n}\right) = 1$ where $\tau(t) = -\log\left(\frac{V(t)}{V(0)}\right)$. For the rescaled solution, we have the rescaled metric, second fundamental form, and curvatures as follows:

- $\widetilde{g_{ij}} = V(t)^{-\frac{2}{n+1}} g_{ij}$ and $\widetilde{h_{ij}} = V(t)^{-\frac{1}{n+1}} h_{ij}$
- $\tilde{H} = V(t)^{\frac{1}{n+1}} H$ and $\tilde{K} = V(t)^{\frac{n}{n+1}} K$
- $\tilde{S} = V(t)^{-\frac{1}{n+1}} S$ and $\tilde{\eta} = V(t)^{\frac{n(\alpha-1)}{n+1}} \eta$ where $\eta(t) = \int_{\mathbb{S}^n} \frac{1}{K^{1-\alpha}}\, d\sigma_{\mathbb{S}^n}$ .

Then we obtain the following corollary.

**Corollary 4.2.** *For the strictly convex and smooth rescaled solution $\tilde{\Sigma}_\tau = \tilde{X}(\mathbb{S}^n, \tau)$ of $\alpha$-Gauss curvature flow, we have the evolution equation of $\tilde{X}$:*

$$\frac{\partial \tilde{X}}{\partial \tau} = -\frac{\tilde{K}^\alpha}{\tilde{\eta}}\tilde{\nu} + \frac{1}{n+1}\tilde{X} \quad \text{on } \mathbb{S}^n \times [0, \infty),$$

*where $\tilde{K}$ is the Gauss curvature and $\tilde{\nu}$ is the unit outward normal of $\tilde{\Sigma}_\tau$.*

*Proof.* Lemma 4.1 implies

$$V(t) = V(0) - \int_0^t \eta(s)\, ds.$$

Since $\frac{\partial \tilde{X}}{\partial \tau} = \frac{\partial \tilde{X}}{\partial t}\frac{dt}{d\tau} = -\frac{K^\alpha V}{\eta V^{1/(n+1)}}\tilde{\nu} + \frac{1}{n+1}\tilde{X}$, we get the result

$$\frac{\partial \tilde{X}}{\partial \tau} = -\frac{\tilde{K}^\alpha}{\tilde{\eta}}\tilde{\nu} + \frac{1}{n+1}\tilde{X}\,. \tag{4.2}$$

□

Now we introduce an integral quantity to analyze the asymptotic behavior of the rescaled hypersurface $\tilde{\Sigma}$.

**Lemma 4.3.** *Let us define the integral quantity $\tilde{\mathcal{I}}$ as follows:*

$$\tilde{\mathcal{I}}(\tau) = \begin{cases} \left( \displaystyle\int_{\mathbb{S}^n} \frac{1}{\tilde{S}^{\frac{1}{\alpha}-1}}\, d\sigma_{\mathbb{S}^n} \right)^{\operatorname{sgn}(\alpha-1)} & \text{for } \alpha > 0 \text{ and } \alpha \neq 1, \\ \displaystyle\int_{\mathbb{S}^n} \log \tilde{S}\, d\sigma_{\mathbb{S}^n} & \text{for } \alpha = 1. \end{cases} \tag{4.3}$$

*Then it satisfies*

$$\frac{\partial}{\partial \tau}\tilde{\mathcal{I}}(\tau) \leq 0. \tag{4.4}$$

*Moreover, the equality holds if and only if $\tilde{K}^\alpha = C\tilde{S}$ a.e. for some positive constant $C$.*

*Proof.*
Case 1. Let us assume that $0 < \alpha < 1$.
By the definition of the rescaled support function $\tilde{S}$ and (4.2), we know that

$$\tilde{K}^{-\alpha}\left( \frac{\partial \tilde{S}}{\partial \tau} - \frac{1}{n+1}\tilde{S} \right) = -\frac{1}{\tilde{\eta}}\,. \tag{4.5}$$

Multiplying both sides of the equation (4.5) by $\tilde{S}^{-\beta}$, where $\beta$ will be chosen later on, implies

$$\frac{1}{\tilde{S}^{\beta}}\left(\frac{\partial \tilde{S}}{\partial \tau} - \frac{1}{n+1}\tilde{S}\right) = -\frac{\tilde{K}^{\alpha}}{\tilde{\eta}\tilde{S}^{\beta}}\,,$$

from the derivation of $\tilde{\mathcal{I}}(\tau)$ with respect to $\tau$, we have

$$\frac{\alpha}{1-\alpha}\left(\mathcal{I}(\tau)\right)^{-2}\frac{\partial}{\partial \tau}\tilde{\mathcal{I}}(\tau) = \int_{\mathbb{S}^n}\frac{\tilde{S}_{\tau}}{\tilde{S}^{\beta}}\,d\sigma_{\mathbb{S}^n} = \frac{1}{n+1}\int_{\mathbb{S}^n}\tilde{S}^{1-\beta}\,d\sigma_{\mathbb{S}^n} - \int_{\mathbb{S}^n}\frac{\tilde{K}^{\alpha}}{\tilde{\eta}\tilde{S}^{\beta}}\,d\sigma_{\mathbb{S}^n} \le 0\,. \tag{4.6}$$

Since $\tilde{\eta}(\tau)$ is positive, (4.6) is non-positive if

$$\frac{1}{n+1}\left(\int_{\mathbb{S}^n}\tilde{S}^{1-\beta}\,d\sigma_{\mathbb{S}^n}\right)\left(\int_{\mathbb{S}^n}\frac{1}{\tilde{K}^{1-\alpha}}\,d\sigma_{\mathbb{S}^n}\right) \le \int_{\mathbb{S}^n}\frac{\tilde{K}^{\alpha}}{\tilde{S}^{\beta}}\,d\sigma_{\mathbb{S}^n}, \tag{4.7}$$

which implies non-positivity of the evolution equation of $\tilde{\mathcal{I}}(\tau)$. Hence it will suffice to show that inequality (4.7) holds. First, notice that we have

$$\begin{aligned}\int_{\mathbb{S}^n}\tilde{S}^{1-\beta}\,d\sigma_{\mathbb{S}^n} &= \int_{\mathbb{S}^n}\left(\tilde{S}^{-\beta}\tilde{K}^{\alpha}\right)^{\frac{\beta-1}{\beta}}\left(\tilde{K}^{-\alpha}\right)^{\frac{\beta-1}{\beta}}\,d\sigma_{\mathbb{S}^n}\\ &\le \left(\int_{\mathbb{S}^n}\tilde{K}^{\alpha}\tilde{S}^{-\beta}\,d\sigma_{\mathbb{S}^n}\right)^{\frac{\beta-1}{\beta}}\left(\int_{\mathbb{S}^n}\frac{1}{\tilde{K}^{\alpha(\beta-1)}}\,d\sigma_{\mathbb{S}^n}\right)^{\frac{1}{\beta}}\end{aligned} \tag{4.8}$$

for $\alpha(\beta-1) = 1-\alpha$, that is, $\beta = \frac{1}{\alpha}$ from the Hölder inequality, which implies

$$\int_{\mathbb{S}^n}\frac{1}{\tilde{S}^{\frac{1}{\alpha}-1}}\,d\sigma_{\mathbb{S}^n} \le \left(\int_{\mathbb{S}^n}\tilde{K}^{\alpha}\tilde{S}^{-\frac{1}{\alpha}}\,d\sigma_{\mathbb{S}^n}\right)^{1-\alpha}\left(\int_{\mathbb{S}^n}\frac{1}{\tilde{K}^{1-\alpha}}\,d\sigma_{\mathbb{S}^n}\right)^{\alpha}. \tag{4.9}$$

We also have

$$\begin{aligned}\int_{\mathbb{S}^n}\frac{1}{\tilde{K}^{1-\alpha}}\,d\sigma_{\mathbb{S}^n} &= \int_{\mathbb{S}^n}\left(\frac{\tilde{S}}{\tilde{K}}\right)^{1-\alpha}\tilde{S}^{-1+\alpha}\,d\sigma_{\mathbb{S}^n}\\ &\le \left(\int_{\mathbb{S}^n}\frac{\tilde{S}}{\tilde{K}}\,d\sigma_{\mathbb{S}^n}\right)^{1-\alpha}\left(\int_{\mathbb{S}^n}\frac{1}{\tilde{S}^{\frac{1}{\alpha}-1}}\,d\sigma_{\mathbb{S}^n}\right)^{\alpha}.\end{aligned} \tag{4.10}$$

Now from (4.8) and (4.10), we get

$$\begin{aligned}&\left(\int_{\mathbb{S}^n}\frac{1}{\tilde{S}^{\frac{1}{\alpha}-1}}\,d\sigma_{\mathbb{S}^n}\right)\left(\int_{\mathbb{S}^n}\frac{1}{\tilde{K}^{1-\alpha}}\,d\sigma_{\mathbb{S}^n}\right)^{1-\alpha}\\ &\qquad\le \left(\int_{\mathbb{S}^n}\tilde{K}^{\alpha}\tilde{S}^{-\frac{1}{\alpha}}\,d\sigma_{\mathbb{S}^n}\right)^{1-\alpha}\left(\int_{\mathbb{S}^n}\frac{1}{\tilde{K}^{1-\alpha}}\,d\sigma_{\mathbb{S}^n}\right)\\ &\qquad\le \left(\int_{\mathbb{S}^n}\tilde{K}^{\alpha}\tilde{S}^{-\frac{1}{\alpha}}\,d\sigma_{\mathbb{S}^n}\right)^{1-\alpha}\left(\int_{\mathbb{S}^n}\frac{\tilde{S}}{\tilde{K}}\,d\sigma_{\mathbb{S}^n}\right)^{1-\alpha}\left(\int_{\mathbb{S}^n}\frac{1}{\tilde{S}^{\frac{1}{\alpha}-1}}\,d\sigma_{\mathbb{S}^n}\right)^{\alpha}\end{aligned}$$

and then

$$\left(\int_{\mathbb{S}^n} \frac{1}{\tilde{S}^{\frac{1}{\alpha}-1}}\, d\sigma_{\mathbb{S}^n}\right)\left(\int_{\mathbb{S}^n} \frac{1}{\tilde{K}^{1-\alpha}}\, d\sigma_{\mathbb{S}^n}\right) \le \left(\int_{\mathbb{S}^n} \tilde{K}^{\alpha}\tilde{S}^{-\frac{1}{\alpha}}\, d\sigma_{\mathbb{S}^n}\right)\left(\int_{\mathbb{S}^n} \frac{\tilde{S}}{\tilde{K}}\, d\sigma_{\mathbb{S}^n}\right) = (n+1)\left(\int_{\mathbb{S}^n} \tilde{K}^{\alpha}\tilde{S}^{-\frac{1}{\alpha}}\, d\sigma_{\mathbb{S}^n}\right) \tag{4.11}$$

since the normalized volume $\tilde{V}(\tau) = \frac{1}{n+1}\left(\int_{\mathbb{S}^n} \frac{\tilde{S}}{\tilde{K}}\, d\sigma_{\mathbb{S}^n}\right) = 1$. The last inequality (4.11) completes the proof of the desired result.

Case 2. Assume $\alpha = 1$.
Since $\tilde{S}$ satisfies the equation $\tilde{S}_\tau = -\frac{\tilde{K}}{|\mathbb{S}^n|} + \frac{1}{n+1}\tilde{S}$ where $|\mathbb{S}^n|$ means the volume of $\mathbb{S}^n$, $\frac{\partial \tilde{\mathcal{I}}(\tau)}{\partial \tau} = \int_{\mathbb{S}^n} \frac{\tilde{S}_\tau}{\tilde{S}}\, d\sigma_{\mathbb{S}^n} \le 0$ is equivalent to

$$\frac{|\mathbb{S}^n|^2}{n+1} \le \int_{\mathbb{S}^n} \frac{\tilde{K}}{\tilde{S}}\, d\sigma_{\mathbb{S}^n}\,. \tag{4.12}$$

Then, we know that

$$|\mathbb{S}^n| \le \left(\int_{\mathbb{S}^n} \frac{\tilde{K}}{\tilde{S}}\, d\sigma_{\mathbb{S}^n}\right)^{\frac{1}{2}}\left(\int_{\mathbb{S}^n} \frac{\tilde{S}}{\tilde{K}}\, d\sigma_{\mathbb{S}^n}\right)^{\frac{1}{2}} = (n+1)^{\frac{1}{2}}\left(\int_{\mathbb{S}^n} \frac{\tilde{K}}{\tilde{S}}\, d\sigma_{\mathbb{S}^n}\right)^{\frac{1}{2}} \tag{4.13}$$

from the Hölder inequality and $\tilde{V}(\tau) = 1$. This implies (4.12) directly. In addition, the equality in (4.4) holds if and only if the equalities hold in (4.9) and (4.13), which implies the equation $\tilde{K}^{\alpha} = C\tilde{S}$ *a.e.* for some positive constant $C$.

□

We can observe that $\tilde{\mathcal{I}}$ is bounded below from [F] for $\alpha = 1$ and $\tilde{\mathcal{I}} \ge 0$ for $\alpha > 0$ and $\alpha \ne 1$. Lemma 4.3 for the evolution equation of $\tilde{\mathcal{I}}$ gives us the following convergence.

**Corollary 4.4.** *For the integral quantity $\tilde{\mathcal{I}}(\tau)$ given by* (4.3), *we have*

$$\lim_{\tau\to\infty} \tilde{\mathcal{I}}(\tau) = \tilde{\mathcal{I}}_0$$

*for some constant $\tilde{\mathcal{I}}_0$, moreover*

$$\lim_{\tau\to\infty} \frac{\partial}{\partial\tau}\tilde{\mathcal{I}}(\tau) = 0.$$

**Lemma 4.5.** *Let us assume that $\tilde{\Sigma}_1$ and $\tilde{\Sigma}_2$ are n-dimensional hypersurfaces embedded in $\mathbb{R}^{n+1}$ and integral quantities of $\tilde{\Sigma}_1$ and $\tilde{\Sigma}_2$ are $\tilde{\mathcal{I}}_1$ and $\tilde{\mathcal{I}}_2$, respectively. If $\tilde{\Sigma}_1 \subset \tilde{\Sigma}_2$, then we have*

$$\tilde{\mathcal{I}}_1 \leq \tilde{\mathcal{I}}_2.$$

*Proof.* Let us $\tilde{S}_1(z)$ and $\tilde{S}_2(z)$ be the support functions of $\tilde{\Sigma}_1$ and $\tilde{\Sigma}_2$, respectively. We know that if $\tilde{\Sigma}_1 \subset \tilde{\Sigma}_2$, then $\tilde{S}_1(z) \leq \tilde{S}_2(z)$. Then by definition of $\tilde{\mathcal{I}}$, we have

$$\begin{aligned}\tilde{\mathcal{I}}_1^{-1} &= \int_{\mathbb{S}^n} \frac{1}{\tilde{S}_1^{\frac{1}{\alpha}-1}} d\sigma_{\mathbb{S}^n} = \int_{\mathbb{S}^n} \frac{1}{\langle z, \tilde{X}_1(\nu^{-1}(z))\rangle^{\frac{1}{\alpha}-1}} d\sigma_{\mathbb{S}^n} \\ &\geq \int_{\mathbb{S}^n} \frac{1}{\langle z, \tilde{X}_2(\nu^{-1}(z))\rangle^{\frac{1}{\alpha}-1}} d\sigma_{\mathbb{S}^n} = \int_{\mathbb{S}^n} \frac{1}{\tilde{S}_2^{\frac{1}{\alpha}-1}} d\sigma_{\mathbb{S}^n} = \tilde{\mathcal{I}}_2^{-1} \quad \text{for } z \in \mathbb{S}^n,\end{aligned}$$

in the case $\alpha > 0$ and $\alpha \neq 1$. Also, for $\alpha = 1$, the desired result follows from $\log \tilde{S}_1(z) \leq \log \tilde{S}_2(z)$.

□

Now we shall show that $\tilde{\Sigma}(\tau)$ has a finite width.

**Lemma 4.6.** *Let us consider an ellipsoid $E(\tau)$ such that $r_{\min}(\tau)$ is equal to half of the minor axis and $r_{\max}(\tau)$ is equal to half of the major axis. Assume that $E(\tau)$ has a fixed volume $V(\tau)$. If $r_{\max}(\tau)$ goes to infinity, then $\tilde{\mathcal{I}}(\tau)$ is also infinity.*

*Proof.* Set $r_1 \cdots r_{n+1} = C$ where $C$ is some positive constant and let $g(x_1, \ldots, x_{n+1}) = \frac{x_1^2}{r_1^2} + \cdots + \frac{x_{n+1}^2}{r_{n+1}^2}$. The equation for the ellipsoid is

$$\frac{x_1^2}{r_1^2} + \cdots + \frac{x_{n+1}^2}{r_{n+1}^2} = 1$$

where $r_1 = r_{\min}$, $r_{n+1} = r_{\max}$ and $r_1 \leq r_2 \leq \cdots \leq r_{n+1}$. Then an ellipsoid can be parameterized by:

$$X = (r_1 q_1, \ldots, r_{n+1} q_{n+1})$$

where $q = (q_1, \ldots, q_{n+1}) \in \mathbb{S}^n$. We also can obtain a normal vector $\boldsymbol{N} = \frac{1}{2}\nabla g = \left(\frac{x_1}{r_1^2}, \ldots, \frac{x_{n+1}}{r_{n+1}^2}\right)$, a unit normal vector $\nu = \frac{\boldsymbol{N}}{\|\boldsymbol{N}\|}$, and the support function $\tilde{S} = \tilde{X} \cdot \tilde{\nu} = \frac{1}{\|\boldsymbol{N}\|}$. Now we have $\frac{x_i}{r_i^2} = \boldsymbol{N}_i = \|\boldsymbol{N}\|\nu_i = \|\boldsymbol{N}\| z_i = \frac{1}{\tilde{S}} z_i$ where $z = (z_1, \ldots, z_n) \in \mathbb{S}^n$, and then

$$\frac{x_i}{r_i} = \frac{r_i}{\tilde{S}} z_i.$$

Since $1 = \frac{x_1^2}{r_1^2} + \cdots + \frac{x_{n+1}^2}{r_{n+1}^2} = \frac{r_1^2}{\tilde{S}^2} z_1^2 + \cdots + \frac{r_{n+1}^2}{\tilde{S}^2} z_{n+1}^2$, we get

$$\tilde{S}^2 = r_1^2 z_1^2 + \cdots + r_{n+1}^2 z_{n+1}^2$$

and we also have

$$\tilde{\mathcal{I}}^{-1} = \int_{\mathbb{S}^n} \frac{1}{\tilde{S}^{\frac{1}{\alpha}-1}} \, d\sigma_{\mathbb{S}^n} = \int_{\mathbb{S}^n} \frac{1}{\left(\sqrt{r_1^2 z_1^2 + \cdots + r_{n+1}^2 z_{n+1}^2}\right)^{\frac{1}{\alpha}-1}} \, d\sigma_{\mathbb{S}^n} .$$

We consider the following case in general: there is $1 \le k \le n+1$ such that $r_{n+1} \ge \cdots \ge r_k \gg r_{k+1} \ge \cdots \ge r_1$ with $r_1 \cdots r_{n+1} = C$, where $C$ is some positive constant. Then we have

$$\begin{aligned} C_1 r_{n+1}^{1-\frac{1}{\alpha}} &\le \int_{\mathbb{S}^n \cap \{\frac{1}{2} \le z_{n+1} \le 1\}} \frac{1}{\left(\sqrt{n+1}\sqrt{r_{n+1}^2 z_{n+1}^2}\right)^{\frac{1}{\alpha}-1}} \, d\sigma_{\mathbb{S}^n} \\ &\le \tilde{\mathcal{I}}^{-1} \le \int_{\mathbb{S}^n} \frac{1}{\left(\sqrt{r_{n+1}^2 z_{n+1}^2}\right)^{\frac{1}{\alpha}-1}} \, d\sigma_{\mathbb{S}^n} \le C_2 r_{n+1}^{1-\frac{1}{\alpha}}, \end{aligned}$$

where $C_1$ and $C_2$ are positive constants. Since $C r_{n+1}^{1-\frac{1}{\alpha}}$ goes to zero for $\alpha < 1$ as $r_{max}(\tau)$ goes to infinity, $\tilde{\mathcal{I}}(\tau)$ is also infinite. Similarly, for $\alpha = 1$, since

$$\int_{\mathbb{S}^n} \Big( \log r_{n+1} + \log |z_{n+1}| \Big) \, d\sigma_{\mathbb{S}^n} \le \tilde{\mathcal{I}} = \int_{\mathbb{S}^n} \log \tilde{S} \, d\sigma_{\mathbb{S}^n} \le c \log r_{n+1}$$

for some positive constant $c$, $\tilde{\mathcal{I}}(\tau) \to \infty$ as $r_{max}(\tau) \to \infty$.

□

Now we shall introduce a theorem called John's Theorem.

**Theorem 4.7** (John's Theorem, [B]). *Let $K$ be a convex body in $\mathbb{R}^n$. Then there exists a unique ellipsoid $E$ of maximal volume which is contained in each $K$. This ellipsoid $E$ is $B_2^n = \left\{x \in \mathbb{R}^n : \sum_1^n x_i^2 \le 1\right\}$ if and only if the following conditions are fulfilled:*

*(1) $B_2^n \subset K$.*

*(2) There are positive numbers $(c_i)_1^m$ and Euclidean unit vectors $(u_i)_1^m$ on the boundary of $K$ such that $\sum_{i=1}^m c_i u_i = 0$ and $\sum_{i=1}^m c_i \langle x, u_i \rangle^2 = |x|^2$ for all $x \in \mathbb{R}^n$.*

We define the width of the convex hypersurface by the function $\tilde{w}(z) = \tilde{S}(z) + \tilde{S}(-z)$ for $z \in \mathbb{S}^n$ and let $\tilde{w}_{\max} = \max_{z \in \mathbb{S}^n} \tilde{w}(z)$ and $\tilde{w}_{\min} = \min_{z \in \mathbb{S}^n} \tilde{w}(z)$. Similarly, set $\tilde{S}_{max} = \max_{z \in \mathbb{S}^n} \tilde{S}(z)$ and $\tilde{S}_{min} = \min_{z \in \mathbb{S}^n} \tilde{S}(z)$. Then we have the following.

**Corollary 4.8.** *For the rescaled hypersurface $\tilde{\Sigma}$ with the normalized volume, there exist some positive constants $0 < c \le C < \infty$ such that*

$$c \le \tilde{w}_{\min} \le \tilde{w}_{\max} \le C$$

*for all $\tau \in [0, \infty)$.*

*Proof.* We know that there exists a unique ellipsoid $E_n$ of maximal volume enclosed by the given convex body $\tilde{\Sigma}$ by Theorem 4.7. Thus we can set up $\tilde{\Sigma}$ between two ellipsoids by using an affine transformation. In other words,

$$E_n \subset \tilde{\Sigma} \subset \sqrt{n}\, E_n\,.$$

Then if the maximum radius of ellipsoid $E_n$ is infinite, the integral quantity $\tilde{\mathcal{I}}$ for $E_n$ is also infinite by Lemma 4.6. This fact and Lemma 4.5 give us that $\tilde{\Sigma}$ does not have the finite integral quantity $\tilde{\mathcal{I}}$. It is a contradiction to Corollary 4.4. Then this implies the desired conclusion for the rescaled hypersurface $\tilde{\Sigma}$ with the normalized volume.

□

**Corollary 4.9.** *For the rescaled hypersurface $\tilde{\Sigma}$ with the normalized volume, we have*

$$\tilde{c} \le \tilde{S}_{min} \le \tilde{S}_{max} \le \tilde{C}$$

*for some constants $0 < \tilde{c} \le \tilde{C} < \infty$ and all $\tau \in [0, \infty)$.*

*Proof.* From Corollary 4.8, we get $\tilde{S}_{max} \le \tilde{C}$ for some positive constant $\tilde{C}$, which implies $\tilde{S}_{min} \ge \tilde{c} > 0$ for some constant $\tilde{c}$ since $\tilde{V}(\tau) = 1$.

□

**Lemma 4.10.**
*If $\tilde{\Sigma}_0$ is strictly convex, then there is a constant $C > 0$ such that*

$$\sup_{z\in\mathbb{S}^n,\ \tau\ge 0} \tilde{K}^{\alpha}(z,\tau) \le C = \max\left( \sup_{z\in\mathbb{S}^n} \tilde{K}^{\alpha}(z,0),\ \Big(\frac{n\alpha+1}{n\alpha\tilde{\rho}_0}\Big)^{n\alpha} \right)$$

*where $\tilde{\rho}_0 = \dfrac{1}{4}\tilde{w}_{\min}$.*

*Proof.* From the evolution equation of $K^{\alpha}$, we have

$$\frac{\partial \tilde{K}^{\alpha}}{\partial \tau} = \frac{\alpha}{\tilde{\eta}}\tilde{K}^{\alpha}\widetilde{(h^{-1})}^{ij}\overline{\nabla}_i\overline{\nabla}_j\tilde{K}^{\alpha} + \frac{\alpha}{\tilde{\eta}}\tilde{K}^{2\alpha}\tilde{H} - \frac{n\alpha}{n+1}\tilde{K}^{\alpha}.$$

By Corollary 4.8, we can consider $\tilde{\rho}_0 = \frac{1}{4}\tilde{w}_{\min}$ and then apply the maximum principle to the function $\tilde{\varphi} = \frac{\tilde{K}^{\alpha}}{\tilde{S}-\rho_0}$. Let us assume that the maximum of $\tilde{\varphi}$ is achieved at the interior point $\tilde{P}_0$ of $\tilde{X}$. Then we have the following properties

$$\tilde{\varphi}_\tau \geq 0,\ \overline{\nabla}_i\tilde{\varphi} = 0 \text{ and } \overline{\nabla}_i\overline{\nabla}_j\tilde{\varphi} \leq 0$$

at $P_0$. Using the evolution equations of $\tilde{K}^\alpha$ and $\tilde{S}$ and calculating by the similar way to Lemma 3.2 implies

$$\begin{aligned} 0 &\leq \frac{\alpha\tilde{K}^{2\alpha}(n-\tilde{S}\tilde{H})}{\tilde{\eta}(\tilde{S}-\tilde{\rho}_0)} + \frac{\alpha}{\tilde{\eta}}\tilde{K}^{2\alpha}\tilde{H} + \frac{\tilde{K}^{2\alpha}}{\tilde{\eta}(\tilde{S}-\tilde{\rho}_0)} - \frac{n\alpha\tilde{K}^\alpha}{n+1} - \frac{\tilde{K}^\alpha\tilde{S}}{(n+1)(\tilde{S}-\tilde{\rho}_0)} \\ &\leq \frac{\tilde{K}^{2\alpha}}{\tilde{\eta}(\tilde{S}-\tilde{\rho}_0)}\Big(n\alpha - \alpha\tilde{\rho}_0\tilde{H} + 1\Big) \end{aligned}$$

at $P_0$, which gives us that

$$0 \leq (n\alpha + 1) - \alpha\tilde{\rho}_0\tilde{H} + 1.$$

Following the same line of the last argument in Lemma 3.2, we get

$$\sup_{z\in\mathbb{S}^n,\ \tau\geq 0} \tilde{K}^\alpha(z,\tau) \leq C = \max\left(\sup_{z\in\mathbb{S}^n}\tilde{K}^\alpha(z,0),\ \Big(\frac{n\alpha+1}{n\alpha\tilde{\rho}_0}\Big)^{n\alpha}\right).$$

□

**Lemma 4.11.** *There is a uniform constant* $0 < \Lambda < \infty$ *such that*

*(i)* $\dfrac{1}{\Lambda} \leq \tilde{S} \leq \Lambda$,

*(ii)* $\dfrac{1}{\Lambda} \leq \tilde{\eta} \leq \Lambda$ *and*

*(iii)* $\dfrac{1}{\Lambda^n} \leq \tilde{K} \leq \Lambda^n$.

*Proof.* (i) $\frac{1}{\Lambda_1} \leq \tilde{S} \leq \Lambda_1$ comes from Corollary 4.9 for some $\Lambda_1 > 0$.

(ii) From Lemma 4.10, we can derive that $\tilde{\eta}(\tau) \geq \frac{1}{\Lambda_l}$ for some positive constant $\Lambda_l > C^{\frac{1-\alpha}{\alpha}}|\mathbb{S}^n|^{-1}$, where $|\mathbb{S}^n|$ is the volume of $\mathbb{S}^n$ and $C$ is the upper bound of $\tilde{K}^\alpha$. In addition, by the Hölder inequality and $\tilde{V} = 1$, we have

$$\tilde{\eta} = \int_{\mathbb{S}^n} \tilde{K}^{\alpha-1}\, d\sigma_{\mathbb{S}^n} \leq \left(\int_{\mathbb{S}^n} \frac{1}{\tilde{K}}\, d\sigma_{\mathbb{S}^n}\right)^{1-\alpha} \cdot |\mathbb{S}^n|^\alpha \leq \Big((n+1)\Lambda_1\Big)^{1-\alpha}|\mathbb{S}^n|^\alpha < \Lambda_u \tag{4.14}$$

for some positive constant $\Lambda_u$. Then we get $\frac{1}{\Lambda_2} \leq \tilde{\eta} \leq \Lambda_2$ by selecting $\Lambda_2 = \max\Big(\Lambda_l,\ \Lambda_u\Big)$.

(iii) Let us consider the evolution of $\overline{S} = \mu\tilde{S}$ for $\mu > 0$. Let $\overline{K}$ and $\overline{H}$ be the Gauss curvature and mean curvature of the hypersurface given by the support function $\overline{S}$, respectively. Then $\overline{K} = \frac{1}{\mu^n}\tilde{K}$, $\overline{H} = \frac{1}{\mu}\tilde{H}$, $\overline{(h^{-1})}^{ij} = \frac{1}{\mu}\widetilde{(h^{-1})}^{ij}$, and $\overline{\eta} = \mu^{(1-\alpha)n}\tilde{\eta}$. Let $\overline{Z}(\tau) = \inf_{z\in\mathbb{S}^n} \overline{K}(z,\tau)$. Then we assume that the interior minimum of $\overline{Z}(\tau)$ is achieved at $\tilde{P}_0 = (z_0, \tau_0)$. From the evolution equation of $\overline{Z}(\tau)$, we have, at $\tilde{P}_0$,

$$\begin{aligned}\frac{\partial\overline{Z}}{\partial\tau} &= \frac{\alpha\mu^{n+1}}{\overline{\eta}}\overline{Z}^{\alpha}\overline{(h^{-1})}^{ij}\overline{\nabla}_i\overline{\nabla}_j\overline{Z} + \frac{\alpha(\alpha-1)\mu^{n+1}}{\overline{\eta}}\overline{Z}^{\alpha-1}\overline{(h^{-1})}^{ij}\overline{\nabla}_i\overline{Z}\,\overline{\nabla}_j\overline{Z} + \frac{\mu^{n+1}}{\overline{\eta}}\overline{Z}\,\overline{H} - \frac{n}{n+1}\overline{Z}\\ &\geq \frac{\mu^{n+1}}{\overline{\eta}}\overline{Z}\,\overline{H} - \frac{n}{n+1}\overline{Z}\\ &\geq \frac{n\mu^{n+1}}{\overline{\eta}}\overline{Z}^{1+1/n} - \frac{n}{n+1}\overline{Z}\\ &\geq n\overline{Z}\Big(\frac{\mu^{n+1}}{\overline{\Lambda}_2}\overline{Z}^{1/n} - \frac{1}{n+1}\Big),\end{aligned}$$

where $\overline{\Lambda}_2 = \mu^{n(1-\alpha)}\Lambda_2$. Set $Q(\tau) = \frac{\mu^{n+1}}{\overline{\Lambda}_2}\overline{Z}^{1/n}(\tau) - \frac{1}{n+1}$ and choose $\mu > \left(\frac{\Lambda_2}{n+1}\right)^{\frac{1}{n\alpha}}\left(\tilde{Z}(0)\right)^{-\frac{1}{n^2\alpha}}$ for $\tilde{Z}(\tau) = \inf_{z\in\mathbb{S}^n}\tilde{K}(z,\tau)$ and $\overline{Z}(\tau) = \frac{1}{\mu^n}\tilde{Z}(\tau)$, which tells us $Q(0) > 0$. Then the evolution equation of $Q(\tau)$ is

$$\begin{aligned}\frac{\partial Q}{\partial\tau} &= \frac{\alpha\mu^{n+1}}{\overline{\eta}}\overline{Z}^{\alpha}\overline{(h^{-1})}^{ij}\overline{\nabla}_i\overline{\nabla}_j Q + \Big(\alpha - \frac{1}{n}\Big)\frac{\alpha n\overline{\Lambda}_2}{\overline{\eta}}\overline{Z}^{\alpha-\frac{1}{n}}\overline{(h^{-1})}^{ij}\overline{\nabla}_i Q\,\overline{\nabla}_j Q\\ &\quad + \frac{\overline{\Lambda}_2^{n-1}}{\mu^{(n+1)(n-1)}}\left(Q + \frac{1}{n+1}\right)^n\left(\frac{\mu^{n+1}}{n\overline{\eta}}\overline{H} - \frac{1}{n+1}\right)\\ &\geq \frac{\overline{\Lambda}_2^{n-1}}{\mu^{(n+1)(n-1)}}\left(Q + \frac{1}{n+1}\right)^n\left(\frac{\mu^{n+1}}{n\overline{\eta}}\overline{H} - \frac{1}{n+1}\right)\\ &\geq \frac{\overline{\Lambda}_2^{n-1}}{\mu^{(n+1)(n-1)}}\left(Q + \frac{1}{n+1}\right)^n Q\end{aligned}$$

at the interior minimum point since $n\overline{Z}^{\frac{1}{n}} \leq n\overline{K}^{\frac{1}{n}} \leq \overline{H}$. By the maximum principle, we have

$$Q(\tau) \geq Q(0) > 0$$

for all $\tau > 0$, which implies $\frac{\partial\overline{Z}}{\partial\tau} > 0$ at $\tilde{P}_0$ and then it gives us contradiction. Hence we obtain

$$\inf_{z\in\mathbb{S}^n}\overline{K}(z,\tau) \geq \inf_{z\in\mathbb{S}^n}\overline{K}(z,0) > 0,$$

and we also have the desired result $\inf_{z\in\mathbb{S}^n}\tilde{K}(z,\tau) \geq \inf_{z\in\mathbb{S}^n}\tilde{K}(z,0) > 0$ for all $\tau$. Combining with Lemma 4.10 implies $\frac{1}{\Lambda_3^n} \leq \tilde{K} \leq \Lambda_3^n$ for some positive constant $\Lambda_3$. Now we select $\Lambda = \max_{i=1,2,3}\Lambda_i$.

□

To obtain the regularity of the solution around the maximal time $T^*$, let us consider the evolution equation (4.2). Then the evolution equation for $\tilde{S}$ is

$$\frac{\partial \tilde{S}}{\partial \tau} = -\frac{\tilde{K}^\alpha}{\tilde{\eta}} + \frac{1}{n+1}\tilde{S}, \tag{4.15}$$

so

$$\Big(\frac{1}{n+1}\tilde{S} - \tilde{S}_\tau\Big)\Big[\det(\overline{\nabla}_i\overline{\nabla}_j\tilde{S} + \tilde{S}\delta_{ij})\Big]^\alpha = \frac{1}{\tilde{\eta}}. \tag{4.16}$$

Now we shall derive $C^{1,1}$-estimates for the solution of (4.16) as in [GH] and [GS].

**Lemma 4.12.** *Suppose that $\tilde{S} \in C^4$ is a solution of the equation* (4.16). *Then we have*

$$\left|\overline{\nabla}^2\tilde{S}\right| \le C \quad on \quad \mathbb{S}^n \times [0,\infty)$$

*where $C$ is a positive constant depending on $\tilde{S}$ and the first derivative of $\tilde{S}$ in time and space.*

*Proof.* Let

$$v(z,\tau) = \left|\tilde{S}(z,\tau)\right|\Big(\overline{\nabla}_\zeta\overline{\nabla}_\zeta\tilde{S}(z,\tau) + \tilde{S}(z,\tau)\Big)\exp\Big(\frac{1}{2}\mu\left|\overline{\nabla}_\zeta\tilde{S}(z,\tau)\right|^2 - \rho\tilde{S}(z,\tau)\Big),$$

where $\mu$ and $\rho$ are positive constants. If the maximum of $v$ is achieved at the initial time, we are done. So we assume that $v$ has its space-time maximum at some interior point $P_0 = (z_0,\tau_0)$ and for some unit vector $\zeta$. We can assume $\zeta = (1,0,\dots,0)$ by choosing an orthonormal frame about $z_0$ and then the matrix $\{\overline{\nabla}_i\overline{\nabla}_j\tilde{S}(z_0)\}$ is diagonal. Then

$$v(z,\tau) = \left|\tilde{S}(z,\tau)\right|\Big(\overline{\nabla}_1\overline{\nabla}_1\tilde{S}(z,\tau) + \tilde{S}(z,\tau)\Big)\exp\Big(\frac{1}{2}\mu\left|\overline{\nabla}_1\tilde{S}(z,\tau)\right|^2 - \rho\tilde{S}(z,\tau)\Big).$$

Let $\mathcal{L}$ be the linearized operator at $P_0$

$$\mathcal{L} = \frac{1}{\Big(\tilde{S}_\tau - \frac{1}{n+1}\tilde{S}\Big)(P_0)}\frac{\partial}{\partial\tau} + \alpha F_{ij}(\overline{\nabla}_k\overline{\nabla}_l\tilde{S} + \tilde{S}\delta_{kl})\overline{\nabla}_i\overline{\nabla}_j.$$

Then $\{\overline{\nabla}_k\overline{\nabla}_l\tilde{S} + \tilde{S}\delta_{kl}\}$ is diagonal. We know that for $F(M) = \log(\det M)$ where $M$ is a positive definite matrix,

$$(F_{ij}) = \frac{\partial F}{\partial M_{ij}} = M^{-1} \quad \text{and} \quad \frac{\partial^2 F}{\partial M_{ij}\partial M_{kl}} = F_{ij,kl} = -F_{ik}F_{jl}.$$

Let

$$\begin{aligned} w &= \log v(z,\tau) \\ &= \log\left|\tilde{S}(z,\tau)\right| + \log\Big(\overline{\nabla}_1\overline{\nabla}_1\tilde{S}(z,\tau) + \tilde{S}(z,\tau)\Big) + \frac{1}{2}\mu\left|\overline{\nabla}_1\tilde{S}(z,\tau)\right|^2 - \rho\tilde{S}(z,\tau) \end{aligned}$$

which also attains its maximum at $P_0$, so $\overline{\nabla} w(P_0)=0$, $\overline{\nabla}_i\overline{\nabla}_j w(P_0)\le 0$, and $w_\tau(P_0)\ge 0$. Since $\overline{\nabla}_i\overline{\nabla}_j\tilde{S}+\tilde{S}\delta_{ij}>0$ from the strict convexity of $\tilde{\Sigma}$, $F_{ij}((\overline{\nabla}_i\overline{\nabla}_j\tilde{S}+\tilde{S}\delta_{ij})(z_0))$ is diagonal, so

$$\mathcal{L}(w)(P_0)=\frac{1}{(\tilde{S}_\tau-\frac{1}{n+1}\tilde{S})(P_0)}\frac{\partial w}{\partial\tau}(P_0)+\alpha F_{ii}((\overline{\nabla}_k\overline{\nabla}_l\tilde{S}+\tilde{S}\delta_{kl})(P_0))\overline{\nabla}_i\overline{\nabla}_j w(P_0)\le 0.$$

From now on, we will use the notation $\overline{\nabla}_{ij}$ in place of $\overline{\nabla}_i\overline{\nabla}_j$ for convenience. We have that at $P_0$

$$\overline{\nabla}_i w=\frac{\overline{\nabla}_i\tilde{S}}{\tilde{S}}+\frac{\overline{\nabla}_{i11}\tilde{S}+\overline{\nabla}_i\tilde{S}}{\overline{\nabla}_{11}\tilde{S}+\tilde{S}}-\rho\overline{\nabla}_i\tilde{S}=0\quad\text{for }\ i=2,\dots,n\,. \tag{4.17}$$

In addition, we get

$$\begin{aligned}\overline{\nabla}_i\overline{\nabla}_i w=&\frac{\overline{\nabla}_{ii}\tilde{S}}{\tilde{S}}-\frac{\left(\overline{\nabla}_i\tilde{S}\right)^2}{\tilde{S}^2}+\frac{\overline{\nabla}_{ii11}\tilde{S}+\overline{\nabla}_{ii}\tilde{S}}{\overline{\nabla}_{11}\tilde{S}+\tilde{S}}-\frac{(\overline{\nabla}_{i11}\tilde{S}+\overline{\nabla}_i\tilde{S})^2}{(\overline{\nabla}_{11}\tilde{S}+\tilde{S})^2}+\mu(\overline{\nabla}_{i1}\tilde{S})^2+\mu\overline{\nabla}_1\tilde{S}\overline{\nabla}_{ii1}\tilde{S}-\rho\overline{\nabla}_{ii}\tilde{S}\\ &\le 0\quad\text{for all }\ i\,,\end{aligned} \tag{4.18}$$

and

$$w_\tau=\frac{\tilde{S}_\tau}{\tilde{S}}+\frac{\overline{\nabla}_{11}\tilde{S}_\tau+\tilde{S}_\tau}{\overline{\nabla}_{11}\tilde{S}+\tilde{S}}+\mu\overline{\nabla}_1\tilde{S}\,\overline{\nabla}_1\tilde{S}_\tau-\rho\tilde{S}_\tau\ge 0 \tag{4.19}$$

at the point $P_0$. Then

$$\begin{aligned}\mathcal{L}(w)(P_0)=&\frac{1}{(\tilde{S}_\tau-\frac{1}{n+1}\tilde{S})(P_0)}\left(\frac{\tilde{S}_\tau}{\tilde{S}}+\frac{\overline{\nabla}_{11}\tilde{S}_\tau+\tilde{S}_\tau}{\overline{\nabla}_{11}\tilde{S}+\tilde{S}}+\mu\overline{\nabla}_1\tilde{S}\,\overline{\nabla}_1\tilde{S}_\tau-\rho\tilde{S}_\tau\right)\\ &+\alpha F_{ii}((\overline{\nabla}_k\overline{\nabla}_l\tilde{S}+\tilde{S}\delta_{kl})(P_0))\Bigg(\frac{\overline{\nabla}_{ii}\tilde{S}}{\tilde{S}}-\frac{\left(\overline{\nabla}_i\tilde{S}\right)^2}{\tilde{S}^2}+\frac{\overline{\nabla}_{ii11}\tilde{S}+\overline{\nabla}_{ii}\tilde{S}}{\overline{\nabla}_{11}\tilde{S}+\tilde{S}}-\frac{(\overline{\nabla}_{i11}\tilde{S}+\overline{\nabla}_i\tilde{S})^2}{(\overline{\nabla}_{11}\tilde{S}+\tilde{S})^2}\\ &+\mu(\overline{\nabla}_{i1}\tilde{S})^2+\mu\overline{\nabla}_1\tilde{S}\overline{\nabla}_{ii1}\tilde{S}-\rho\overline{\nabla}_{ii}\tilde{S}\Bigg)\le 0.\end{aligned} \tag{4.20}$$

Since $\left(\frac{1}{n+1}\tilde{S}-\tilde{S}_\tau\right)\left[\det(\overline{\nabla}_i\overline{\nabla}_j\tilde{S}+\tilde{S}\delta_{ij})\right]^\alpha=\frac{1}{\tilde{\eta}}$, after differentiation and then some calculations, we have

$$\frac{\frac{1}{n+1}\overline{\nabla}_1\tilde{S}-\overline{\nabla}_1\tilde{S}_\tau}{\frac{1}{n+1}\tilde{S}-\tilde{S}_\tau}+\alpha F_{ii}\left\{\overline{\nabla}_{1ii}\tilde{S}+\overline{\nabla}_1\tilde{S}\right\}=0\quad\text{at }\ P_0. \tag{4.21}$$

Once again the differentiation implies

(4.22)
$$\frac{\frac{1}{n+1}\overline{\nabla}_{11}\tilde{S}-\overline{\nabla}_{11}\tilde{S}_\tau}{\frac{1}{n+1}\tilde{S}-\tilde{S}_\tau}-\frac{(\frac{1}{n+1}\overline{\nabla}_1\tilde{S}-\overline{\nabla}_1\tilde{S}_\tau)^2}{(\frac{1}{n+1}\tilde{S}-\tilde{S}_\tau)^2}+\frac{\alpha(\overline{\nabla}_{11ii}\tilde{S}+\overline{\nabla}_{11}\tilde{S})}{\overline{\nabla}_{ii}\tilde{S}+\tilde{S}}-\frac{\alpha(\overline{\nabla}_{1ij}\tilde{S}+\overline{\nabla}_1\tilde{S}\,\delta_{ij})^2}{(\overline{\nabla}_{ii}\tilde{S}+\tilde{S})(\overline{\nabla}_{jj}\tilde{S}+\tilde{S})}=0\quad\text{at } P_0.$$

We use the properties of covariant derivatives:

$$\overline{\nabla}_{kji}\tilde{S}=\overline{\nabla}_{jik}\tilde{S}+\delta_{ik}\overline{\nabla}_j\tilde{S}-\delta_{ij}\overline{\nabla}_k\tilde{S} \tag{4.23}$$

and

$$\overline{\nabla}_{lkji}\tilde{S}=\overline{\nabla}_{jilk}\tilde{S}+2\delta_{kl}\overline{\nabla}_{ji}\tilde{S}-2\delta_{ij}\overline{\nabla}_{lk}\tilde{S}+\delta_{il}\overline{\nabla}_{jk}\tilde{S}-\delta_{kj}\overline{\nabla}_{li}\tilde{S}. \tag{4.24}$$

After using the formulas (4.21)-(4.24) and the following properties

$$\frac{\alpha(\overline{\nabla}_{ij1}\tilde{S}+\delta_{j1}\overline{\nabla}_i\tilde{S})^2}{(\overline{\nabla}_{ii}\tilde{S}+\tilde{S})(\overline{\nabla}_{jj}\tilde{S}+\tilde{S})}=\frac{\alpha(\overline{\nabla}_{i11}\tilde{S}+\overline{\nabla}_i\tilde{S})^2}{(\overline{\nabla}_{ii}\tilde{S}+\tilde{S})(\overline{\nabla}_{11}\tilde{S}+\tilde{S})}+\sum_{i=1}^n\sum_{j=2}^n\frac{\alpha\left(\overline{\nabla}_{ij1}\tilde{S}\right)^2}{(\overline{\nabla}_{ii}\tilde{S}+\tilde{S})(\overline{\nabla}_{jj}\tilde{S}+\tilde{S})}$$

and

$$\frac{\alpha\mu(\overline{\nabla}_{11}\tilde{S}+\tilde{S})\overline{\nabla}_1\tilde{S}\,\delta_{i1}\overline{\nabla}_i\tilde{S}}{\overline{\nabla}_{ii}\tilde{S}+\tilde{S}}=\alpha\mu(\overline{\nabla}_1\tilde{S})^2$$

and several computations, (4.20) can be simplified to

(4.25)
$$\begin{aligned}0\geq&-\frac{\tilde{S}_\tau(\overline{\nabla}_{11}\tilde{S}+\tilde{S})}{\tilde{S}\left(\frac{1}{n+1}\tilde{S}-\tilde{S}_\tau\right)}-\frac{\tilde{S}_\tau+\overline{\nabla}_{11}\tilde{S}}{\frac{1}{n+1}\tilde{S}-\tilde{S}_\tau}+\frac{\left(\frac{1}{n+1}\overline{\nabla}_1\tilde{S}-\overline{\nabla}_1\tilde{S}_\tau\right)^2}{\left(\frac{1}{n+1}\tilde{S}-\tilde{S}_\tau\right)^2}+\frac{\alpha(\overline{\nabla}_{11}\tilde{S}-\overline{\nabla}_{ii}\tilde{S})}{\overline{\nabla}_{ii}\tilde{S}+\tilde{S}}\\&+\sum_{i=1}^n\sum_{j=2}^n\frac{\alpha\left(\overline{\nabla}_{ij1}\tilde{S}\right)^2}{(\overline{\nabla}_{ii}\tilde{S}+\tilde{S})(\overline{\nabla}_{jj}\tilde{S}+\tilde{S})}+\frac{\alpha\overline{\nabla}_{ii}\tilde{S}(\overline{\nabla}_{11}\tilde{S}+\tilde{S})}{\tilde{S}(\overline{\nabla}_{ii}\tilde{S}+\tilde{S})}-\frac{\alpha\left(\overline{\nabla}_i\tilde{S}\right)^2(\overline{\nabla}_{11}\tilde{S}+\tilde{S})}{\tilde{S}^2(\overline{\nabla}_{ii}\tilde{S}+\tilde{S})}-\alpha\mu(\overline{\nabla}_1\tilde{S})^2\\&-\frac{\frac{\mu}{n+1}(\overline{\nabla}_{11}\tilde{S}+\tilde{S})(\overline{\nabla}_1\tilde{S})^2}{\frac{1}{n+1}\tilde{S}-\tilde{S}_\tau}+\frac{\rho(\overline{\nabla}_{11}\tilde{S}+\tilde{S})\tilde{S}_\tau}{\frac{1}{n+1}\tilde{S}-\tilde{S}_\tau}+\frac{\alpha\mu(\overline{\nabla}_{11}\tilde{S}+\tilde{S})(\overline{\nabla}_{i1}\tilde{S})^2}{\overline{\nabla}_{ii}\tilde{S}+\tilde{S}}-\frac{\alpha\rho(\overline{\nabla}_{11}\tilde{S}+\tilde{S})\overline{\nabla}_{ii}\tilde{S}}{\overline{\nabla}_{ii}\tilde{S}+\tilde{S}}.\end{aligned}$$

In addition, since

(4.26)
$$\begin{aligned}&\sum_{i=1}^n\sum_{j=2}^n\frac{\alpha\left(\overline{\nabla}_{ij1}\tilde{S}\right)^2}{(\overline{\nabla}_{ii}\tilde{S}+\tilde{S})(\overline{\nabla}_{jj}\tilde{S}+\tilde{S})}-\frac{\alpha\left(\overline{\nabla}_i\tilde{S}\right)^2(\overline{\nabla}_{11}\tilde{S}+\tilde{S})}{\tilde{S}^2(\overline{\nabla}_{ii}\tilde{S}+\tilde{S})}\\=&\sum_{i=2}^n\sum_{j=2}^n\frac{\alpha\left(\overline{\nabla}_{ij1}\tilde{S}\right)^2}{(\overline{\nabla}_{ii}\tilde{S}+\tilde{S})(\overline{\nabla}_{jj}\tilde{S}+\tilde{S})}+\sum_{i=2}^n\frac{2\alpha\rho\left(\overline{\nabla}_{1i1}\tilde{S}\right)\overline{\nabla}_i\tilde{S}}{\overline{\nabla}_{ii}\tilde{S}+\tilde{S}}-\sum_{i=2}^n\frac{\alpha\rho^2(\overline{\nabla}_{11}\tilde{S}+\tilde{S})\left(\overline{\nabla}_i\tilde{S}\right)^2}{\overline{\nabla}_{ii}\tilde{S}+\tilde{S}}-\frac{\alpha\left(\overline{\nabla}_1\tilde{S}\right)^2}{\tilde{S}^2}\end{aligned}$$

from (4.23) and (4.17), we have

(4.27)

$$0 \geq -\frac{\tilde{S}_\tau(\overline{\nabla}_{11}\tilde{S}+\tilde{S})}{\tilde{S}\left(\frac{1}{n+1}\tilde{S}-\tilde{S}_\tau\right)} - \frac{\tilde{S}_\tau+\overline{\nabla}_{11}\tilde{S}}{\frac{1}{n+1}\tilde{S}-\tilde{S}_\tau} + \frac{\left(\frac{1}{n+1}\overline{\nabla}_1\tilde{S}-\overline{\nabla}_1\tilde{S}_\tau\right)^2}{\left(\frac{1}{n+1}\tilde{S}-\tilde{S}_\tau\right)^2} + \frac{\alpha\overline{\nabla}_{11}\tilde{S}}{\tilde{S}} + \sum_{i=2}^{n}\sum_{j=2}^{n}\frac{\alpha\left(\overline{\nabla}_{ij1}\tilde{S}\right)^2}{(\overline{\nabla}_{ii}\tilde{S}+\tilde{S})(\overline{\nabla}_{jj}\tilde{S}+\tilde{S})}$$
$$+\sum_{i=2}^{n}\frac{2\alpha\rho\left(\overline{\nabla}_{1i1}\tilde{S}\right)\overline{\nabla}_i\tilde{S}}{\overline{\nabla}_{ii}\tilde{S}+\tilde{S}} - \sum_{i=2}^{n}\frac{\alpha\rho^2(\overline{\nabla}_{11}\tilde{S}+\tilde{S})\left(\overline{\nabla}_i\tilde{S}\right)^2}{\overline{\nabla}_{ii}\tilde{S}+\tilde{S}} - \frac{\alpha\left(\overline{\nabla}_1\tilde{S}\right)^2}{\tilde{S}^2} - \alpha\mu(\overline{\nabla}_1\tilde{S})^2$$
$$-\frac{\frac{\mu}{n+1}(\overline{\nabla}_{11}\tilde{S}+\tilde{S})(\overline{\nabla}_1\tilde{S})^2}{\frac{1}{n+1}\tilde{S}-\tilde{S}_\tau} + \frac{\rho(\overline{\nabla}_{11}\tilde{S}+\tilde{S})\tilde{S}_\tau}{\frac{1}{n+1}\tilde{S}-\tilde{S}_\tau} + \alpha\mu\left(\overline{\nabla}_{11}\tilde{S}\right)^2 - \frac{\alpha\rho(\overline{\nabla}_{11}\tilde{S}+\tilde{S})\overline{\nabla}_{ii}\tilde{S}}{\overline{\nabla}_{ii}\tilde{S}+\tilde{S}}.$$

Let $\gamma_i = \overline{\nabla}_{ii}\tilde{S}+\tilde{S}$. Then (4.27) can be written as

(4.28)

$$0 \geq -\alpha - \frac{\tilde{S}_\tau\gamma_1}{\tilde{S}\left(\frac{1}{n+1}\tilde{S}-\tilde{S}_\tau\right)} - \frac{\gamma_1}{\frac{1}{n+1}\tilde{S}-\tilde{S}_\tau} + \frac{\alpha\gamma_1}{\tilde{S}} + \sum_{i=2}^{n}\frac{\alpha\rho\left(\overline{\nabla}_i\tilde{S}\right)^2\gamma_1}{\gamma_i}\left(\rho-\frac{2}{\tilde{S}}\right) - \frac{\alpha\left(\overline{\nabla}_1\tilde{S}\right)^2}{\tilde{S}^2}$$
$$-\alpha\mu(\overline{\nabla}_1\tilde{S})^2 - \frac{\frac{\mu}{n+1}(\overline{\nabla}_1\tilde{S})^2\gamma_1}{\frac{1}{n+1}\tilde{S}-\tilde{S}_\tau} + \frac{\rho\tilde{S}_\tau\gamma_1}{\frac{1}{n+1}\tilde{S}-\tilde{S}_\tau} + \alpha\mu{\gamma_1}^2 - 2\alpha\mu\tilde{S}\gamma_1 + \alpha\mu\tilde{S}^2 - \alpha\rho\gamma_1 + \frac{\alpha\rho\tilde{S}\gamma_1}{\gamma_i}$$

at $P_0$. We obtained the lower and upper bounds of $\tilde{S}$ on $[0,\infty)$ in Lemma 4.11, and $|\overline{\nabla}_i\tilde{S}|$ is also bounded for $i = 1,\dots,n$ since $\tilde{\Sigma}$ is strictly convex. In addition, since $\frac{1}{n+1}\tilde{S}-\tilde{S}_\tau$ has the positive lower bound from Lemma 4.11, choosing $\mu$ and $\rho$ such that

$$0 \geq A{\gamma_1}^2 + B\gamma_1 + C_1,$$

where $A$ is a positive constant and $B$ and $C_1$ are some constants, give us the desired result.

□

**Corollary 4.13.** *There exists some positive constant C such that*

$$\sup_{x\in\tilde{\Sigma},\tau\geq 0} \widetilde{\mathcal{H}} \leq C.$$

*Moreover, $\tilde{\lambda}_{min} \geq C_1 > 0$ for some constant $C_1$. Here $\tilde{\lambda}_{min} = \tilde{\lambda}_1 \leq \cdots \leq \tilde{\lambda}_n = \tilde{\lambda}_{max}$ where $\tilde{\lambda}_i'$s are the eigenvalues of $(\tilde{h}^i_j)$.*

Furthermore, combining Lemma 4.10 and Corollary 4.13 implies the following Corollary.

**Corollary 4.14.** *All curvatures on the rescaled hypersuface* $\tilde{\Sigma}$ *are bounded above and below by the uniform constants. In other words, there exists some constant* $0 < M < \infty$ *such that*

$$\frac{1}{M} \leq \tilde{\lambda}_{min} \leq \tilde{\lambda}_{max} \leq M.$$

## 5. Existence of Solutions and Proof of Main Theorem

**5.1. Short time existence.** Let us assume that $\Sigma_t$ is smooth. Then we get the uniform $C^{1,1}$-estimates of the coefficients of our equation (2.4) and this equation becomes uniformly parabolic. Thus the regularity theory of uniform parabolic equations and application of the implicit function theorem give us the short time existence as in [Li].

**5.2. Long time existence.** Let $\lambda_i$ be the eigenvalues of $(h^i_j)$. We know that $\lambda_i$ is positive by the strict convexity. Also, we have $K = \lambda_1 \cdots \lambda_n \leq C_1$ and $\mathcal{H} = \frac{1}{\lambda_1} + \cdots + \frac{1}{\lambda_n} \leq C_2$ from Lemma 3.2 and Lemma 3.3, where $\lambda_1 \leq \lambda_2 \leq \cdots \leq \lambda_n$ and $C_1$ and $C_2$ are some positive constants. These give us, for each $i = 1, \ldots, n$,

$$0 < \frac{1}{C_2} \leq \lambda_i$$

from $\frac{1}{\lambda_i} < \frac{1}{\lambda_1} + \cdots + \frac{1}{\lambda_n} \leq C_2$ and also

$$0 < \lambda_i \leq \frac{C_1}{\Pi_{j \neq i} \lambda_j} \leq C_1 C_2^{n-1},$$

which imply there are $0 < \lambda \leq \Lambda < \infty$ satisfying

$$\lambda |\xi|^2 \leq K^{\alpha} (h^{-1})^{ij} \xi_i \xi_j \leq \Lambda |\xi|^2.$$

Then we know that the support function $S(z, t)$ satisfies a uniformly parabolic equation in $\Sigma_t$. Hence $S(z, t)$ is $C^{2,\gamma}$ and then $C^{\infty}$ in $\Sigma_t$ through the standard bootstrap argument using the Schauder theory. If there is a $0 < T_1 < T^*$ such that $\Sigma_t$ is smooth on $[0, T_1)$ but not smooth after $T_1$, the uniform $C^{2,\gamma}$-estimates for $S(z, t)$ implies that $\Sigma_{T_1}$ is $C^{2,\gamma}$, and therefore $C^{\infty}$. From the short time existence and uniqueness, $\Sigma_t$ is $C^{\infty}$ on $[0, T_1 + \delta)$ for some small $\delta > 0$. It is a contradiction. Therefore $T_1 = T^*$ and there is a smooth solution $\Sigma_t$ on $[0, T^*)$. Also, the solution $\Sigma_t$ will be strictly convex by Lemma 3.1.

**Proof of Theorem 1.1.** We have the uniform bounds of curvature and all of the higher derivatives of the second fundamental form to the rescaled manifold by Corollary 4.14 and then the equation (4.16) will be uniformly parabolic. In addition, we have $C^{1,1}$-regularity of the solution $\tilde{S}$ from Lemma 4.12. By applying the Harnack inequality to the linearized equation satisfied by $\tilde{S}_\tau$, we obtain that $\tilde{S}_\tau$ is Hölder continuous through a similar argument as in [GH]. We can apply the Evans-Krylov theorem and the Schauder estimates (see [CC]) to the concave operator obtained by taking exponent $\frac{1}{n\alpha}$ to the equation (4.16), which implies $C^{2,\gamma}$-regularity of $\tilde{S}$ for $0 < \gamma < 1$. And then we have the smooth and strictly convex rescaled solution by the standard bootstrap argument using the Schauder theory and Corollary 4.13. In other words, for every sequence of $\tau_k \to \infty$, we can find a subsequence $\tau_{k_i}$ such that $\tilde{S}(\cdot, \tau_{k_i}) \to \tilde{S}_*(\cdot)$. Also, the integral quantity

$$\tilde{\mathcal{I}}(\tau) = \begin{cases} \left( \int_{\mathbb{S}^n} \frac{1}{\tilde{S}^{\frac{1}{\alpha}-1}} \, d\sigma_{\mathbb{S}^n} \right)^{\operatorname{sgn}(\alpha-1)} & \text{for } \alpha > 0 \text{ and } \alpha \neq 1, \\ \int_{\mathbb{S}^n} \log \tilde{S} \, d\sigma_{\mathbb{S}^n} & \text{for } \alpha = 1 \end{cases}$$

satisfies the monotonicity $\frac{d}{d\tau}\tilde{\mathcal{I}}(\tau) \leq 0$, and equality holds if and only if $\tilde{K}^\alpha = C\tilde{S}$, for some positive constant $C$, holds for a choice of origin. For the limit manifold $\tilde{\Sigma}_{T^*}$ of the volume rescaled manifold $\tilde{\Sigma}_{\tau_{i_k}}$, following the same argument as in Theorem 16, [A2], $\tilde{\mathcal{I}}(\tau) \to -\infty$ if $\tilde{\Sigma}_{T^*}$ does not satisfy $\tilde{K}_*^\alpha = \tilde{C}_*\tilde{S}_*$ a.e. for some positive constant $\tilde{C}_*$, which gives a contradiction. Therefore the proof is complete.

□

**Acknowledgement.** Ki-Ahm Lee was supported by the National Research Foundation of Korea(NRF) grant funded by the Korea government(MSIP) (No.2014R1A2A2A01004618). Ki-Ahm Lee also hold a joint appointment with the Research Institute of Mathematics of Seoul National University.

## References

[A1] B. Andrews, *Gauss Curvature Flow: The Fate of the Rolling Stones*, Invent. Math. **138** (1999), 151–161.

[A2] B. Andrews, *Motion of hypersurfaces by Gauss curvature*, Pacific J. Math. **195** (2000), 1–34.

[A3] B. Andrews, *Harnack inequalities for evolving hypersurfaces.*, Math. Z., **217** (1994), 179–197.

[A5] B. Andrews, *Contraction of convex hypersurfaces by their affine normal.*, J. Differential Geom., **43** (1996), 207–230.

[B] K. Ball, *An Elementary Introduction to Modern Convex Geometry*, Flavors of geometry, **58**, Math. Sci. Res. Inst. Publ., 31, Cambridge Univ. Press, Cambridge, 1997.

[CC] L.A. Caffarelli, X. Cabré, *Fully nonlinear elliptic equations*, American Mathematical Society, Providence, RI, 1995.

[C1] B. Chow, *Deforming convex hypersurfaces by the nth root of the Gaussian curvature*, J. Differential Geom., **22** (1985), 117–138.

[C2] B. Chow, *Deforming convex hypersurfaces by the square root of the scalar curvature*, Invent. math., **87** (1987), 63–82.

[F] W. Firey, *Shapes of worn stones*, Mathematica **21** (1974), 1–11.

[GH] C. Gutiérrez, Q. Huang, *A generalization of a theorem by Calabi to the parabolic Monge-Ampère equation*, Indiana Univ. Math. Journal, **47** (1998), 1459–1480.

[GaH] M. Gage, R. S. Hamilton, *The heat equation shrinking convex plane curves*, J. Differential Geom., **23** (1986), 69–96.

[GS] B. Guan, J. Spruck, *Boundary-value problems on $\mathbb{S}^n$ for surfaces of constant Gauss curvature*, Ann. of Math., **138** (1993), 601–624.

[KLR] L. Kim, K. Lee, E. Rhee, *$\alpha$-Gauss Curvature flows with flat sides*, J. of Differential Equations, **254** (2013), 1172–1192.

[KS] N.V. Krylov, N.V. Safonov, *Certain properties of solutions of parabolic equations with measurable coefficients*, Izvestia Akad. Nauk. SSSR, **40** (1980), 161–175.

[Li] G.M. Lieberman, *Second Order Parabolic Differential Equations*, World Scientific Publishing Co., Inc., River Edge, NJ, 1996.

[T] K. Tso, *Deforming a hypersurface by its Gauss-Kronecker curvature*, Comm. Pure appl. Math. **38** (1985), 867–882.
[Z] X.-P. Zhu, *Lectures on Mean Curvature Flows*, Studies in Advanced Mathematics, **32**, American Mathematical Society, Providence, RI, 2002.

Department of Mathematical Sciences, Seoul National University, Seoul 151-747, Republic of Korea & Department of Mathematics, Hokkaido University, Sapporo 060-0810, Japan
*E-mail address*: lmkim.math@gmail.com

Department of Mathematical Sciences, Seoul National University, Seoul 151-747, Republic of Korea & School of Mathematics, Korea Institute for Advanced Study, Seoul 130-722, Republic of Korea
*E-mail address*: kiahm@snu.ac.kr